\documentclass[a4paper, 12pt]{article}
\usepackage{amsmath, amsfonts, amsthm, amssymb, qtree, enumitem, hyperref}

\theoremstyle{definition} 
\theoremstyle{plain} 
\theoremstyle{definition} 
\theoremstyle{plain} 
\theoremstyle{plain} 
\theoremstyle{plain} 

\DeclareMathOperator{\wneg}{\sim}

\begin{document}
\title{G\"odel's Theorem and Direct Self-Reference}
\author{Saul A. Kripke}
\date{}

\maketitle

\begin{abstract}
In his paper on the incompleteness theorems, G\"odel seemed to say that a direct way of constructing a formula that says of itself that it is unprovable might involve a faulty circularity. In this note, it is proved that `direct' self-reference can actually be used to prove his result. 

\noindent \textit{Keywords/MSC2010 Classification Codes}: 03-03, 03B10, direct self-reference, G\"odel numbering, Incompleteness.
\end{abstract}

In his epochal 1931 paper on the incompleteness theorems G\"odel wrote, regarding his undecidable statement, saying, in effect, ``I am unprovable,''
\begin{quote}
Contrary to appearances, such a proposition involves no faulty circularity, for initially it [only] asserts that a certain well-defined formula (namely, the one obtained from the $q$th formula in the lexicographic order by a certain substitution) is unprovable. Only subsequently (and so to speak by chance) does it turn out that this formula is precisely the one by which the proposition itself was expressed.\cite[p.~151, fn.~15]{Godel 1931}
\end{quote}
The tendency of this footnote might be that if one tried to express the statement directly, there might indeed be a faulty circularity. In my own paper, ``Outline of a Theory of Truth'', I expressed some doubt about this tendency, and stated that the G\"odel theorem could be proved using `direct' self-reference. To quote this paper,
\begin{quote}
A simpler, and more direct, form of self-reference uses demonstratives or proper names: Let `Jack' be a name of the sentence `Jack is short', and we have a sentence that says of itself that it is short. I can see nothing wrong with ``direct'' self-reference of this type. If `Jack' is not already a name in the language,$^{4}$ why can we not introduce it as a name of any entity we please? In particular, why can it not be a name of the (uninterpreted) finite sequence of marks `Jack is short'? (Would it be permissible to call this sequence of marks ``Harry,'' but not ``Jack''? Surely prohibitions on naming are arbitrary here.) There is no vicious circle in our procedure, since we need not \textit{interpret} the sequence of marks `Jack is short' before we name it. Yet if we name it ``Jack,'' it at once becomes meaningful and true. (Note that I am speaking of self-referential sentences, not self-referential propositions.$^{5}$)

In a longer version, I would buttress the conclusion of the preceding paragraph not only by a more detailed philosophical exposition, but also by a mathematical demonstration that the simple kind of self-reference exemplified by the ``Jack is short'' example could actually be used to prove the G\"odel incompleteness theorem itself (and also, the G\"odel-Tarski theorem on the undefinability of truth). Such a presentation of the proof of the G\"odel theorem might be more perspicuous to the beginner than is the usual one. It also dispels the impression that G\"odel was forced to replace direct self-reference by a more circumlocutory device. The argument must be omitted from this outline.$^{6}$\cite[pp. 77--78]{Kripke 1975}.

\noindent $^{4}$ We assume that `is short' \textit{is} already in the language.

\noindent $^{5}$ It is \textit{not} obviously possible to apply this technique to obtain ``directly'' self-referential \textit{propositions}.

\noindent $^{6}$ There are several ways of doing it, using either a nonstandard G\"odel numbering where statements can contain numerals designating their own G\"odel numbers, or a standard G\"odel numbering, plus added constants of the type of `Jack'.
\end{quote}
As I say in footnote~6, there are several ways of obtaining the G\"odel theorem using direct self-reference, analogously to `Jack is short'. Note that, as I say in footnote~5, I am not claiming that this technique could be used to obtain `direct' self referential \textit{propositions}. Perhaps this is what G\"odel had in mind in his own footnote. Let us look at how G\"odel's first incompleteness theorem is to be done these ways. 

One way, observed independently by Raymond Smullyan, is to use a nonstandard G\"odel numbering where a formula can contain a numeral designating its own G\"odel number.\footnote{Almost everything I say about Smullyan is based on conversation with him. We haven't been able to find a published reference.} This is clearly impossible under the usual G\"odel prime power numbering, or various other variants. That this should not happen might be a natural restriction, since a G\"odel number might be thought to correspond with a formula as a composite object, and this might be thought to preclude a formula from containing a numeral designating its own G\"odel number.

Nevertheless we propose a nonstandard G\"odel numbering allowing a statement to contain a numeral designating its own G\"odel number. We can proceed as follows. Let $x_{1}$ be a fixed variable. Let $A_{1}(x_{1})$, $A_{2}(x_{1})$, \ldots be an enumeration of all those formulae that contain at most $x_{1}$ free (G\"odel's class signs). We assume that the language of the system studied contains, either directly or by virtue of interpretation, the language of arithmetic. Numerals can be assumed to be terms $0$ followed by (or preceded by) finitely many successor symbols (allowing none). We use $0^{(n)}$ for $0$ with $n$ successor symbols; $n=0$ is allowed, as I said. $0^{(n)}$ denotes the number $n$. 

Let the `original' G\"odel numbering be G\"odel's own prime power product numbering, except that the smallest prime used is $3$, so that G\"odel numbers are always odd. In the `new' numbering, all G\"odel numbers coincide with the `original', except that for each $n$, the formula
\begin{equation*}
\label{eqn:2}
(\exists x_{1})(x_{1}=0^{(2k_{n})}\wedge A_{n}(x_{1}))
\end{equation*}
gets the G\"odel number $2k_{n}$, where $k_{n}$ is the original G\"odel number of $(\exists x_{1})(x_{1}=0^{(n)}\wedge A_{n}(x_{1}))$. The `new' numbering allows a formula to contain a numeral designating its G\"odel numbering, and in that sense it is a self-referential G\"odel numbering. 

In this self-referential G\"odel numbering, every formula $A_{n}(x_{1})$ has an `instance' $(\exists x_{1})(x_{1}=0^{(2k_{n})}\wedge A_{n}(x_{1}))$ asserting that its own G\"odel number satisfies $A_{n}(x_{1})$. The G\"odel incompleteness theorem is the special case where $A_{n}(x_{1})$ is unprovability in the system. 

The `new' G\"odel number is in effective $1-1$ correspondence with the old one, where the inverse function is also effective. The complexity of a property in the arithmetical hierarchy ($\Sigma_{1}^{0}, \Pi_{1}^{0}$, etc.) does not change when the `new' G\"odel numbers replace the old.\footnote{The existential quantifier in $(\exists x_{1})(x_{1}=0^{(2k_{n})}\wedge A_{n}(x_{1}))$ does not increase complexity, since the formula is equivalent to $A(0^{(2k_{n})})$. For the G\"odel theorem, the important feature is that `unprovability in the system' is still a $\Pi_{1}^{0}$ predicate. If the language contains non-arithmetical predicates that have hierarchical status, say, in the hyperarithmetical or analytical hierarchy, these also remain the same. 

Further, since the domain and range of both of these $1-1$ mappings are recursive, the mapping can be extended to a recursive permutation of the natural numbers.}

Remark: We could have simply used $A_{n}(0^{(2k_{n})})$. But using $(\exists x_{1})(x_{1}=0^{(2k_{n})}\wedge A_{n}(x_{1}))$ has two advantages, both emphasized by Raymond Smullyan. First, concatenation is easier to arithmetize than substitution, which has some value even here. Second, this choice guarantees that no formula gets two G\"odel numbers, which is not clear if we use the simpler version. It also may not be harmful if some formulae get two or more G\"odel numbers, but in that case we cannot think of formulae as \textit{identified} with their G\"odel numbers. 

Much more natural, in my opinion, for getting `direct' self-reference, is the use of added constants.\footnote{One of the referees has commented, regarding numerals, ``it is highly plausible to view them as proper names; but they have also been seen as analogues of quotational names and structural-descriptive names. Numerals---unlike constants---are complex closed terms and at least in this respect similar to a term built from a term with function symbols for the primitive recursive functions for substitution and numerals''. The doubts in question do not apply to constants, as the referee says. Nevertheless, the use of a self-referential G\"odel numbering, though somewhat artificial, does seem more directly self-referential than the traditional G\"odelian argument.} Once again, let $A_{1}(x_{1}), A_{2}(x_{1}), \ldots, A_{i}(x_{1})$ be an enumeration of all those formulae of the language $L$ containing at most the variable $x_{1}$ free. Let $a_{1}, a_{2}, \ldots$ be a denumerable list of constants, none of which are in the original language $L$. Give G\"odel numbers to the formulae of the extended language in a conventional way. Now add to the axioms of the original system $S$, an infinite set of axioms $a_{i}=0^{(n_{i})}$, where $n_{i}$ is the G\"odel number of $A_{i}(a_{i})$. The resulting system is $S'$. Pretty clearly $S'$ is a conservative extension of $S$. It simply extends $S$ by adding an infinite set of constants with specific numerical values---new names for particular numbers. Every proof in $S'$ becomes a proof in $S$ if each constant is replaced by the corresponding numeral. 

Now, for the G\"odel (first incompleteness) theorem, consider the predicate (in the original language, without the extra constants) that says that $x_{1}$ is the G\"odel number of a formula unprovable in the extended system $S'$. Call this $\wneg\text{Thm}_{S'}(x_{1})$. This is a $\Pi_{1}^{0}$ predicate and is a formula $A_{i}(x_{1})$ for some $i$. In our construction, there is a constant $a_{i}$ of $S'$ such that $a_{i}=0^{(n_{i})}$ is an axiom of $S'$, where $n_{i}$ is the G\"odel number of $A_{i}(a_{i})$. Then it is easily shown, in the usual way, that $A_{i}(a_{i})$ is unprovable in $S'$ if $S'$ is consistent, and that $\wneg A_{i}(a_{i})$ is unprovable in $S'$ if $S'$ is $\omega$-consistent (or even $1$-consistent). Given that $S'$ is a conservative extension of $S$, the same things must be true of $A_{i}(0^{(n_{i})})$ in the original system $S$ as well. 

The particular case of unprovability is just one example. For any definable property we obviously can use this construction to formulate a statement that `says of itself' that it has that property (the `self-reference lemma'). 

A more delicate construction uses only a single constant $a$, and a single extra axiom, in defining the system $S^{*}$ extending $S$. It will have only one extra axiom, which will be called (\ref{eqn:1}). In the case of G\"odel's first incompleteness theorem, we wish to construct a statement saying that it itself is unprovable. So given the constant $a$, we consider the single formula standing for the statement
\begin{equation}
\tag{$*$}
\label{eqn:1}
\parbox{10cm}{$a$ is unprovable in the system $S^{*}$, which is $S$ extended by adding $\ulcorner a=0^{(a)}\urcorner$ as a single axiom}
\end{equation}
(Of course this is to be written out symbolically.) The constant $a$ occurs in (\ref{eqn:1}) only as explicitly shown. Note that in (\ref{eqn:1}), the constant $a$ is both used and mentioned. In a standard G\"odel numbering of the extended language, (\ref{eqn:1}) has a particular G\"odel number. Say it is $n$. Then let the extended system $S^{*}$ be axiomatized by adding $a=0^{(n)}$ as an axiom. 

Given this extra axiom, $a=0^{(n)}$, note that (\ref{eqn:1}) will say that it itself (i.e., its G\"odel number) is unprovable in the system $S^{*}$ itself. Now we can once again deduce that (\ref{eqn:1}) is unprovable in $S^{*}$ if $S^{*}$ is consistent, its negation is unprovable in $S^{*}$ if $S^{*}$ is $\omega$-consistent, etc. Also once again $S^{*}$ is a conservative extension of $S$, so that any proof in $S$ is valid in $S^{*}$ iff it was valid in $S$, just as in the previous argument with an infinity of constants. And once again, this allows us to transfer the result to the original system $S$. Clearly also the argument extends to results about any definable $A(x_{1})$ in place of unprovability (the self-reference lemma). 

\textit{Remark}. In spite of my own use of the Quinean terminology of `use' and `mention' in the preceding paragraph, I am not following his practices of quotation and quasi-quotation. Many formal expressions in the paper are used autonymously. The corners in  (\ref{eqn:1}) are used to denote G\"odel numbers of formulae $\ulcorner a=0^{(a)}\urcorner$, where $a$ is a constant on the left, but on the right it is an as yet unspecified numeral (the constant $a$ is used, rather than mentioned). In other words, `$0^{(a)}$' is not a term in the language, since it means $0$ followed by $a$ successors, and therefore what formula is mentioned in the corners is yet to be specified. In the whole formula in corners, $a$ is `mentioned'.\footnote{A referee has suggested that the reader would profit from the following further readings on direct self-reference: \cite{Heck 2004,Visser 2004,Grabmayr Visser 2020,Halbach Visser 2014a,Halbach Visser 2014b,Picollo 2018}.}

\paragraph{Acknowledgments} My thanks to Yale Weiss for editorial support. Special thanks to Romina Padr\'o, Panu Raatikainen, and Albert Visser for their comments and for encouraging me to write this paper. This paper has been completed with support from the Saul Kripke Center, at The City University of New York, Graduate Center. 

\paragraph{Affiliation} Distinguished Professor in the Philosophy and Computer Science Programs at The Graduate Center, CUNY. \textit{Address} The Saul Kripke Center, The Graduate Center, CUNY, 365 Fifth Ave., Room 7118, New York, NY 10016, USA. \textit{Email} \href{mailto:kripkecenter@gc.cuny.edu}{kripkecenter@gc.cuny.edu}.

\end{document}